\documentclass[11pt]{article}
\usepackage{graphicx}
\usepackage{amssymb}
\usepackage{amsmath}
\usepackage{enumerate}

\newcommand{\N}{\mathbb{N}}

\newtheorem{definition}{{\bf Definition}}[section]
\newtheorem{theorem}[definition]{{\bf Theorem}}

\newtheorem{corollary}[definition]{{\bf Corollary}}
\newtheorem{proposition}[definition]{\noindent {\bf Proposition}}

 \newtheorem{problem}[definition]{\noindent {\bf Problem}}

\newtheorem{remark}[definition]{\noindent {\bf Remark}}
\newtheorem{fact}[definition]{{\bf Fact}}

\newtheorem{lemma}[definition]{\noindent {\bf Lemma}}
\def\endproof{\hfill {\kern 6pt\penalty 500
\raise -0pt\hbox{\vrule \vbox to5pt {\hrule width 5pt
\vfill\hrule}\vrule}}}
\def\centerpicture #1 by #2 (#3){\leavevmode
        \vbox to #2{
        \hrule width #1 height 0pt depth 0pt
        \vfill
        \special{pictfile #3}}}
\begin{document}

\setlength{\parindent}{10 pt}
\title{From well-quasi-ordered sets to better-quasi-ordered sets }
\author{Maurice Pouzet 
\thanks{Supported by Intas}\\
PCS, Universit\'e Claude-Bernard Lyon1,
 Domaine de Gerland -b\^at.\\ Recherche [B], \\
 50 avenue Tony-Garnier, F$69365$ Lyon cedex 07, France\\
{\tt pouzet@univ-lyon1.fr }
\and
Norbert Sauer  \thanks{Supported by NSERC of
Canada Grant \# 691325} \\
 Department of Mathematics and Statistics, The University of Calgary,\\
 Calgary, T2N1N4, Alberta, Canada\\
{\tt nsauer@math.ucalgary.ca}
}
\date{}
\maketitle
\begin{abstract} We consider conditions which force a well-quasi-ordered poset (wqo) to be better-quasi-ordered (bqo). In particular we obtain that if a poset $P$ is wqo and the set $S_{\omega}(P)$ of strictly increasing sequences of  elements of $P$ is bqo under domination,   then $P$ is bqo. As a consequence, we get the same conclusion if  $S_{\omega} (P)$ is replaced by $\mathcal J^1(P)$, the collection of  
non-principal ideals of $P$, or  by  $AM(P)$, the collection of maximal antichains of $P$ ordered by domination.  It then  follows  that an interval order which is wqo is in fact bqo.
\end{abstract}


\section{Introduction and presentation of the results}

\subsection{How to read this paper}

Section \ref{index} contains  a collection of definitions, notations and basic facts. The specialist reader should be able to read the paper with only occasional use of Section \ref{index} to check up on some notation. Section \ref{index} provides readers which are not very familiar with the topic of the paper with  some background,  definitions and simple derivations from those definitions. Such  readers will have to  peruse Section \ref{index} frequently. 

The paper is organized as follows. Section \ref {secbarriers} provides the  basics  behind the notion of bqo posets and  develops the technical tools we need to work with barriers and concludes with the  proof of a result about $\alpha$-bqo's from which Theorem~\ref{nonprincipal} follows.   We present some topological  properties of ideals in Section \ref{sec:ideals} and discuss minimal type posets in Section \ref{sec:min}. The  proof of Theorem \ref{criticbqo} is contained  in Section \ref{sec:pred}.  In Section \ref{maximal} we present  constructions involving maximal antichains of prescribed size.

\subsection{Background}

Since their introduction by G.Higman \cite {higm},  well-quasi-ordered, (wqo), sets have played an important role in several areas of mathematics: algebra (embeddability of free algebras in skew-fields, elimination orderings) set theory and logic (comparison of chains, termination of rewriting systems, decision problems)  analysis (asymptotic computations, symbolic dynamic ). A recent example is given by the Robertson-Seymour Theorem \cite{Rob-Sey} asserting that the collection of finite graphs is well-quasi-ordered by the minor relation.

In this paper we deal with the stronger notion of   better-quasi-ordered, (bqo), posets. Bqo posets where  introduced by C.~St.~J.~A.~Nash-Williams, see \cite{nwbqots},  to prove that the class of infinite trees is  wqo under topological embedding.

Better-quasi-orders  enjoy several properties of well-quasi-orders. 
For example,  finite posets are bqo.  Well ordered chains are bqo, finite unions and finite products of bqo posets are bqo. The property of being   bqo is  preserved under restrictions 
and epimorphic images.  Still there is  a substantial difference: Well-quasi-ordered posets are not preserved under the infinitary construction described in the next paragraph, but better-quasi-ordered posets are .  

A  basic result due to G.Higman, see \cite{higm},   asserts that a poset $P$ is wqo if and only if  ${\bf I}(P)$, the set of initial segments of $P$,  is well-founded. On the other hand,    Rado  \cite {rado} has produced an example of a well-founded partial order  $P$ for which  ${\bf I}(P)$ is  well-founded and contains infinite antichains.   The idea behind the bqo notion is to forbid this  situation: ${\bf I}(P)$ and all its iterates, ${\bf I}({\bf I}(\cdots ({\bf I}(P)\cdots ))$ up to the ordinal $\omega_1$, have to be well-founded and hence wqo.

This  idea  is quite natural but not workable.  (Proving that a two element set satisfies this property is far from being an easy task). The  working definition, based upon the notion of {\it barrier}, invented by  C.~St.~J.~A.~Nash-Williams,  is quite involved, see \cite{nwwqoft} and \cite{nwbqots}.  Even using this working condition,  it is not so  easy to see wether a wqo is a bqo or not. We aim to arrive at a better understanding of bqo posets and   consider two special problems to see if indeed we obtained such a better understanding. 

We solved the first problem, to characterize bqo interval orders, completely, see Theorem \ref{thm:int}. The second was Bonnet's problem, see Problem~\ref{prob:B}. We related the property of a poset to be bqo to the bqo of various  posets associated to a given poset, in particular the poset of the maximal antichains under the domination order.  We think that those results stand on their own but unfortunately don't seem to be strong enough to solve Bonnet's problem.

\subsection {The results}

Let $P$ be a poset.

For $X,Y\subseteq P$ let $X\leq_{dom}Y$ if for every $x\in X$ there is a $y\in Y$ with $x\leq y$. The quasi-order $\leq_{dome}$ is the {\em domination} order on $P$ and $S_\omega(P)$ is the set of strictly increasing $\omega$-sequences of $P$.  We will prove, see Theorem~\ref{nonprincipal2} and the paragraph before Theorem~\ref{nonprincipal2}:

\begin{theorem}\label{nonprincipal} If $P$ is wqo, and $(S_{\omega}(P);\leq_{dom})$ is bqo then  
$P$ is  bqo. 
\end{theorem}

Let $C\in S_\omega(P)$. Then $\downarrow\hskip-2pt C$ is an ideal of $P$. On the other hand if $I$ is an ideal with denumerable cofinality then $I=\downarrow\hskip-2pt C$ for some $C\in S_\omega(P)$. 

Let $\mathcal J^{\neg \downarrow\hskip -2pt }(P)$ be the set of non principal ideals. Since ideals with denumerable cofinality are non-principal, we obtain from Theorem \ref{nonprincipal} and the property of bqo to be preserved under restrictions that:

\begin{corollary}\label{cornonprincipal} If $P$ is wqo and $\mathcal J^{\neg \downarrow\hskip -2pt }(P)$ is bqo then  $P$ is  bqo. 
\end{corollary}

The poset $(S_{\omega}(P);\leq_{dom})$ is often more simple than the poset $P$. So for example
if $P$ is finite $S_\omega(P)=\emptyset$. It follows trivially from Definition \ref{defin:bqo} that the empty poset is bqo and hence from Theorem \ref{nonprincipal} that finite posets are bqo. A result which is of course well known. Also:

\begin{corollary}\label {cor2}If $P$ is wqo and $\mathcal J^{\neg \downarrow\hskip -2pt }(P)$ is finite then $P$ is bqo.
\end{corollary}
Corollary \ref{cornonprincipal}  was conjectured by the first author in his thesis\cite{pouzettr} and a proof of  Corollary~\ref{cor2} given there. The  proof is  given in \cite{fraissetr} Chapter $7$, subsections $7.7.7$ and $7.7.8$. pp $217-219$.  

The above considerations suggest that $S_\omega(P)$ corresponds to some sort of derivative. As already observed,  the elements of  $S_\omega(P)$ generate  the non-principal ideals of $P$ with denumerable cofinality. The set of ideals $\mathcal{J}(P)$ of $P$ form the base set of a Cantor space $C$, see Section \ref{sec:ideals}. It follows that $P$ is finite if and  only if the first Cantor-Bendixson derivative of $C$ is finite and that $\mathcal{J}(P)$ is finite if and only if the second Cantor-Bendixson derivative of $C$ is finite.

The space $C$ contains just one limit if and only if  $\mathcal{J}(P)$ is a singleton space.  Such posets are called minimal type posets. Minimal type posets occur  naturally in symbolic dynamics. See section \ref{sec:min} for details.

\begin{lemma} If $P$ is an interval order then $\mathcal J^{\neg \downarrow\hskip -2pt }(P)$ is a chain. \end{lemma}
\begin{proof} Let $I,J\in \mathcal J^{\neg \downarrow\hskip -2pt }(P)$. If $I\setminus J\not = \emptyset$ and  $J\setminus I \not = \emptyset$ pick $x \in I \setminus J$ and $y\in J\setminus I$. Since $I$ is not a principal ideal then $x$ is not a maximal element in $I$, so we may pick $x'\in I$ such that $ x<x'$. For the same reason, we may pick $y'\in J$ such that $y<y'$. Clearly, the poset induced on $\{x,x',y,y'\}$ is a $\underline 2\oplus \underline  2$. But then $P$ is not an interval order.
\end{proof}

Corollary   \ref {cornonprincipal} has an other immediate  consequence: 
\begin{corollary}\label{chain}
If $P$ is wqo and $\mathcal J^{\neg \downarrow\hskip -2pt } (P)$ is a chain then $P$ is bqo.
\end{corollary}

Indeed, if $P$ is wqo then ${\bf I}(P)$ is well-founded. In particular $\mathcal J^{\neg \downarrow\hskip -2pt }(P)$ is well-founded.  If $\mathcal J^{\neg \downarrow\hskip -2pt }(P)$ is a chain, this is   a well-ordered chain, hence a bqo. From Corollary  \ref{cornonprincipal} $P$ is bqo.

From Corollary \ref{chain}, this gives:
\begin{theorem} \label{thm:int}
An interval order is bqo iff it is wqo. 
\end{theorem}

The following Theorem is an immediate consequence of Theorem \ref{thmsucc}. (See Section \ref{index} for a definition of the notions used in Theorem \ref{criticbqo}.)
\begin{theorem} \label{criticbqo}
Let $P$ be a poset. If $P$ has no infinite antichain, then the following properties are equivalent:
\begin{enumerate}  [{(i)}]
\item $P$ is bqo.
\item $(P;\leq_{succ})$ is bqo.
\item $(P;\leq_{pred})$ is bqo.
\item $(P;\leq_{crit})$ is bqo. 
\item $AM(P)$ is bqo.\label{item:5}
\end{enumerate}
\end{theorem}

As indicated earlier part of the  motivation for this research was an intriguing problem due to Bonnet, see \cite{bonnet}.

\begin{problem}\label{prob:B}
Is every wqo poset  a countable union of bqo posets?
\end{problem}

Item {\em \ref{item:5}} of Theorem \ref{criticbqo}  may suggest to attack Bonnet's problem using the antichains of the poset.    Note that if a poset $P$  is  wqo but  not bqo, then it  contains antichains of arbitrarily large finite size. Indeed, if the size of antichains  of a poset $P$ is bounded by some integer, say $m$, then from Dilworth's theorem, $P$ is the union of at most $m$  chains. If $P$ is well founded, these chains are well ordered chains, hence are bqo, and $P$ is bqo as a finite union of bqo's.

 For  each integer $m$,  let  $AM_{m}(P)$ be the collection of maximal antichains having size $m$ and $\bigcup {AM_{m}(P)}$ be the union of these maximal antichains. Trivially, $AM(P)$  is the union of the sets $ {AM_{m}(P)}$  for $m\in \N$  if and only if $P$ contains no infinite antichain.  It its tempting to use this decomposition to attack Bonnet's problem.  That is,  consider the union  $\bigcup  {AM_{m}(P)}=P$  of a well-quasi-ordered poset $P$. 
 
This does not work in general:  there  are  wqo posets $P$ for which $\bigcup AM_2(P)$ is not bqo ( Lemmas \ref{construct1} and \ref {counterexample}). Still, this works for wqo posets $P$ for which $\bigcup  {AM_{m}(P)}$ is bqo for every $m$ and $AM(P)$ is not bqo (these $P$ are not bqo). Rado's poset provides an example, see Lemma \ref{AM(R)}.
\vskip 4pt 
Looking at  the relationship between $AM_m(P)$ and $\bigcup {AM_m(P)}$, we prove, see Theorem \ref{thm:222}:
\begin{theorem}
Let $P$ be a poset with no infinite antichain,  then $AM_2(P)$ is bqo if and only if $\bigcup {AM_2(P)}$ is bqo. 
\end{theorem}

This does not extend:  it follows from Corollary \ref{construct2} that there exists a wqo poset $P$ for which $AM_{3}(P)$ is  bqo but $\bigcup {AM_{3}(P)}$ is not bqo.

\section{Barriers and better-quasi-orders}\label{secbarriers}

\subsection{Basics} We use Nash-William's notion of  bqo, see \cite{nwwqoft}, and refer to  Milner's exposition of bqo theory, see \cite{milnerwqobqo}. See Section \ref{index} for the basic definitions.

The following result due to F.Galvin(1968) extends the partition theorem of F.P.Ramsey. 
\begin{theorem}\cite{galvin} For every subset $B$ of $[\N ]^{<\omega}$ there is an infinite subset $X$ of $\N$ such that either $[X]^{<\omega}\cap B=\emptyset $ or $[X]^{<\omega}\cap B$ is a block.
\end{theorem}

Trivially, every block contains a thin block,    the set $\min_{\leq _{in}}(B)$ of $\leq_{in}$  minimal elements of the block.  Moreover, if  $B$ is a  block, resp.a thin block,  and $X$ is an infinite subset of  $\bigcup B$ then $B_{\restriction X}$  is a block, resp a thin block.
The  theorem of Galvin implies  the following result of Nash-Williams, see \cite{milnerwqobqo}. 
\begin{theorem} \label{partnash}
\begin{enumerate} [{(a)}]
\item Every block contains a barrier. 
\item For every partition of a barrier into finitely many  parts, one contains a barrier. 
\end{enumerate}
\end{theorem}

 The partial order  $(B,\leq_{lex})$  is the lexicographic sum of the partial orders $(B_{(i)},\leq_{lex})$: 
 $$(B,\leq_{lex})=\sum_{i\in \N}(B_{(i)}, \leq_{lex}).$$

 Let $\mathrm{T}(B)$ be the tree $\mathrm{T}(B):=\bigl(\{t : \exists s\in B(t\leq_{in}s)\}, \leq_{in}\bigr)$ with root $\emptyset$ and $\mathrm{T}^d(B)$ the dual order of $\mathrm{T}(B)$.  If $\mathrm{T}(B)$ does not contain an infinite chain then  $\mathrm{T}^{d}(B)$ is well founded and the height function satisfies
\begin{align*}
& h\bigl(\emptyset ,  \mathrm{T}^{d}(B)\bigr)=   \\
&  \sup\{h\bigl((a),  \mathrm{T}^{d}(B)\bigr) +1: (a)\in \mathrm{T}(B)\} 
    = \sup\{h\bigl(\emptyset,  \mathrm{T}^d(_{(a)}B)\bigr) +1: (a)\in \mathrm{T}(B) \}. 
\end{align*}
Induction on the height gives then   that  $\mathrm{T}(B)$ is   well ordered under the lexicographic order.  The order type of $\mathrm{T}$ being at most $\omega^{\alpha}$ where $\alpha:=h(\emptyset ,  \mathrm{T}^{d}(B))$.   From this fact, we deduce:

 \begin{lemma}\cite{pouz72} Every thin block,  and in particular every barrier,  is well ordered under the lexicographic order. 
\end{lemma}
This allows to associate with every  barrier its order-type.  We note  that $\omega$ is the least possible order-type.  An ordinal  $\gamma$ is the  order-type  of a barrier if and only if  $\gamma=\omega^{\alpha}\cdot n$ where $n<\omega$ and $n=1$ if $\alpha <\omega$ \cite{assous74}. Every barrier contains a barrier whose order-type is an indecomposable ordinal.

\begin{definition}\label{defin:bqo}
A map $f$ from a barrier $B$ into a poset $P$ is {\em good} if there are $s,t\in B$ with $s\lhd t$ and  $f(s) \leq f(t)$. Otherwise $f$ is {\em bad}. \\
Let $\alpha$ be a denumerable ordinal. A poset $P$ is {\em $\alpha$-better-quasi-ordered}  if every map $f: B\rightarrow P$, where $B$ is a barrier of order type at most $\alpha$,  is good.  \\  
A poset $P$ is {\em better-quasi-ordered} if it is $\alpha$-better-quasi-ordered   for every denumerable ordinal $\alpha$. 
\end{definition}

It is known and easy to see that a poset $P$ is $\omega$-better-quasi-ordered   if and only if it is well-quasi-ordered. Remember that  we  abbreviate better-quasi-order  by bqo.
Since every barrier contains a barrier with indecomposable order type, only barriers with indecomposable order type need to be taken into account in the definition of bqo.  In particular, we only need to consider $\alpha$-bqo for indecomposable ordinals $\alpha$.  Note that for indecomposable ordinals  $\alpha$ the notion of $\alpha$-bqo leeds to different objets \cite{marcone1994} .

We will need the following results of Nash-Williams  (for proofs in the context of $\alpha$-bqo, see  \cite{milnerwqobqo} or \cite{pouzet92}):
\begin{lemma}\label{higbqo} Let $P$ and $Q$ be partial orders, then:
\begin{enumerate}[{(a)}]  
\item Finite partial orders and well ordered chains are bqo.
\item   If  $P,Q$ are $\alpha$-bqo then the direct sum $P\oplus Q$ and the direct product $P\times Q$ are $\alpha$-bqo.
\item If $P$ is $\alpha$-bqo and $f:P\longrightarrow Q$ is order-preserving then $f(P)$ is $\alpha$- bqo.
\item If $P$ embeds into $Q$ and $Q$ is $\alpha$-bqo then $P$ is $\alpha$-bqo.
\item If  $C\subseteq \mathfrak {P}(P)$  is $\alpha$-bqo then the set of finite unions of members of $C$ is $\alpha$-bqo.\end{enumerate}
\end{lemma}

It follows from Item $e$ that if $P$ is  $\alpha$-bqo then 
$I_{<\omega}(P)$  is  $\alpha$-bqo which in turn implies that if $P$ is $\alpha$-bqo  then $AM(P)$ is $\alpha$-bqo. (If $P$ is $\alpha$-bqo then it is $\omega$-bqo and hence   well-quasi-ordered and hence does not contain infinite antichains. Then $AM(P)$ embeds into $A(P)$ which in turn embeds into $I_{<\omega}(P)$.)  It follows  from Item $d$ that if $P$ is bqo then every restriction of $P$ to a subset of its elements is also bqo.

\subsection{Barrier constructions}
 Let $B$ be a subset of $[\N]^{<\omega}$.  See Section \ref{index} for notation.  
 
If $B$ is a block then $B^2$ is a block and if $B$ is a thin block then $B^2$ is  a thin block.  Moreover, if $B$ is a thin block, and $u\in B^2$, then there is a unique pair $s,t\in B$ such that $s\lhd t$ and $u= s\cup t$.  If  $B$ is a block, then $\bigcup {_{*}B}=\bigcup{_{*}B}= \bigcup {B}\setminus \min(\bigcup B)$ and  $_{*}B$ is a block. Moreover, if $B$ is well ordered under the lexicographic order then $_{*}B$ is well  ordered too and if the type of $B$ is an indecomposable ordinal $\omega^{\gamma}$
 then the type of   $_{*}B$ is at most $\omega^{\gamma}$.

If  $C$ is a block and $B:= C^2$ then $_{*}B= C\setminus C_{(a)}$, where $a$ is the least element of $\bigcup C$.

The following Lemma is well known and follows easily from the definition.

\begin{lemma}
If $B$ is a barrier, then $B^{2}$ is a barrier and if $B$ has type $\alpha$ then $B^2$ has type $\alpha\cdot\omega$.
\end{lemma}

A generalization, Lemma \ref{lem:Marcone} below,  was given by 
 A.Marcone. We recall his   construction and result, see  \cite {marcone2001} Lemma 8 pp. 343.
 
 Let $B$ be a subset of $[\N]^{<\omega}$. Then $B^{\circ}$ is the set of all elements $s\in B$ with the property that for all $i\in \bigcup{B}$ with $i<s(0)$ there is an element $t\in B$ with $(i)\cdot {_\ast s}\leq_{in} t$. In other words $s\in B^{\circ}$ if $(i)\cdot {_\ast s}\in \mathrm{T}(B)$ for all $i\in \bigcup{B}$ with $i<s(0)$. Let $B^\prime:=\{{_\ast s} : s\in B^\circ\}\setminus \{\emptyset\}$.

\begin{lemma}\label{lem:Marcone} Let $B$ be  a thin block of type larger than $\omega$, then: 
\begin{enumerate}
\item $B'$ is a thin block.
\item  For every  $u\in (B')^2$  there is some $s\in B$ such that $s\leq_{in} u$.
\item If the type of $B$ is at most $\omega^{\gamma}$ then $B'$ contains a barrier  of type at most  $\omega^{\gamma}$ if $\gamma$ is a limit ordinal and at most $\omega^{\gamma-1}$ otherwise. 
\end{enumerate}
\end{lemma}

\begin{remark} We may note that for every $s\in B$ such that $s\leq_{in} u:= s'\cup t'$ we have $s'\leq_{in}s$. Indeed, otherwise $s\leq_{in }s'$, but $s':=_{*}s''$ for some $s''\in B$, hence $s\subseteq s''$ contradicting the fact that $B$ is a barrier.
\end{remark}
A barrier $B$ is {\it end-closed} if 
\begin{equation} \label{eqendclosed}
s\cdot(a)\in B \; \text{ for every  }s\in B_{*} \text{ and } a\in \bigcup B \; \text{ with } a>\lambda (s).
\end{equation}
 For example, $[\N]^{n}$ is end-closed for every $n$,  $n\geq 1$,  as well as the barrier  $B:= \{s\in \N^{<\omega}: l(s)=s(0)+2\}$.

If $B, C$ are two barriers with the same domain, the set $B*C:= \{s\cdot t: s\in C, t\in B, \lambda(s)<t(0)\}$ is a barrier, the  {\it product } of $B$ and $C$, see \cite{pouz72}. Its  order-type is $\omega^{\gamma+\beta}$ if $\omega^{\gamma}$ and $\omega^{\beta}$ are the order-types  of $B$ and $C$ respectively.  For example, the product  $[\bigcup B]^1*B$  is  end-closed.  Provided that $B$ has type $\omega^{\beta}$, it has type $\omega^{1+\beta}$. The converse holds, namely:

\begin{fact}\label{fact0}
The set $D\subseteq \N^{<\omega}$ is an end-closed barrier of type larger than $\omega$ if and only if   $D_{*}$ is a barrier and $D:= [\bigcup {D_{*}}]^1*{D_{*}}$. 
\end{fact}

\begin{lemma}  Every barrier   $B$ contains an end-closed subbarrier $B^\prime$. \end{lemma} 
\begin{proof}
Induction on the order-type $\beta$ of $B$. 

If $\beta:=\omega$ then $B=[\bigcup B]^{1}$ and we may set $B':= B$.

Suppose $\beta>\omega$ and  every barrier of type smaller than $\beta$ contains an end-closed subbarrier. 

The set  $S(B):= \{i \in \bigcup B: (i) \in B\}$  is an initial segment of $\bigcup B$.  (Indeed, let $i\in S(B)$ and  $j<i$ with $j\in \bigcup B$. Select $X\in  [\bigcup B]^{\omega}$ such that $(j,i)\leq_{in}X$. Since $B$ is a barrier, $X$ has an initial segment $s\in B$.  Since $B$ is an antichain w.r.t. inclusion $i\not\in s$, hence  $s=(j)$.) The type of $B$ is larger than $\omega$,  hence $S(B)\not=\bigcup{B}$.  Set $i_{0}:=\min( \bigcup {B}\setminus S(B))$.

The set $_{(i_{0})}B$ is a barrier because $i_{0}\not \in S(B)$.  Hence induction applies providing  some   $X_{0} \subseteq   \bigcup {B}\setminus (S(B)\cup \{i_{0}\})$ such that $\{s:  (i_0)\cdot s \in B\cap [X_{0}]^{<\omega} \}$ is an end-closed barrier of domain $X_{0}$. It follows that $\{i_0\}\cup X_0\subseteq \bigcup{B}$.

Starting with $(x_0,X_0)$ we  construct  a sequence  $( i_{n}, X_{n})_{n<\omega}$ such that for every  $n<\omega$: 
\begin{enumerate}
\item $\{s:  (i_n)\cdot s \in B\cap [X_{n}]^{<\omega} \}$ is an end-closed barrier of domain $X_{n}$.
\item $\{i_{n+1}\}\cup X_{n+1}\subseteq X_{n}$. 
\end{enumerate}
Let  $n<\omega$.  If $(i_m, X_{m})_{m<n}$ is defined for all  $m<n$ replace $B$ by $B\cap [X_{n-1}]^{<\omega}$ in the  construction of $x_0$ and $X_0$ to obtain  $i_{n}$ and $X_{n}$. 

Then for   $X:= \{i_{n}:  n<\omega\}$  set $B':= B\cap [X]^{<\omega}$. 
\end{proof}

 \subsection{ On the comparison of blocks }

Let $B,B'$ be two subsets of $[\N]^{<\omega}$. We write $B'\leq_{in} B$ if :
 \begin{equation}\label{eqbarrier}
\mbox{For every} \;  s'\in B'\;  \mbox { there is some}\;  s\in B\; \mbox { such that } \; s'\leq_{in} s. 
\end{equation}
This is the quasi-order of domination associated with the order $\leq_{in}$ on $[\N]^{<\omega}$. 
 \begin{fact}\label{fact1}
Let $B$, $B'$ be two thin blocks.  If $B'\leq_{in} B$ then for all $s,t\in B$ and $s^\prime, t^\prime\in B^\prime$:
\begin{enumerate}[{(a)}]
\item $\bigcup B'\subseteq \bigcup B$.
\item  If  $s\leq_{in} s'$ then $s=s'$, hence $\leq_{in}$ is a partial order on  thin blocks.
\item If $\bigcup B=\bigcup B'$ then for every $s\in B$ there is some $s'\in B'$ such that $s'\leq_{in} s$. 
\item If $s'\lhd t' $  then $s\lhd t$ for some $s,t\in B$ with $s'\leq_{in} s$ and $t'\leq_{in} t$.
 \item The set $B'':= B'\cup D$ with $D:= \{s\in B: \forall  s'\in B'\,  (s'\not \leq_{in} s)\}$ is a thin block  and $\bigcup B''= \bigcup B$ and $B''\leq_{in}B$.

\end{enumerate}
\end{fact}
\begin{proof} $(a), (b), (c)$ follow from the definitions. 

$(d)$. Let $s', t'\in B'$ with $s'\lhd t' $ and let $t'':= s'\cup t'$. Since  $\bigcup{B'}\subseteq \bigcup B$, $t''\subseteq  \bigcup B$. Let  $X\in [\bigcup B] ^{\omega}$ such that $t''\leq_{in} X$. There are  $s, t\in B$ such that  $s\leq _{in}X$ and $t\leq_{in}{_*}X$. We have $s\lhd t$.   It follows from $(b)$ that $s'\leq_{in} s$ and $t'\leq_{in} t$. 

$(e)$  $\bigcup B''= \bigcup B'\cup \bigcup D$ and  $\bigcup B'\subseteq \bigcup B$ imply  $\bigcup B''\subseteq \bigcup B$.  For the converse, let  $x\in \bigcup B\setminus \bigcup B'$. Since $B$ is a block there is some $s\in B$ having $x$ as first element. Clearly  $s\in D$, hence $x\in D$, proving $\bigcup B''= \bigcup B$. From the definition,  $B''$ is an antichain. Now, let $X\subseteq B$. We prove that some initial segment $s''$ belongs to  $B''$. Since $B$ is a block, some initial segment $s$  of $X$ belongs to $B$. If $s\in D$ set $s'':= s$. Otherwise some initial segment $s'$ of $s$ is in $B'$. Set $s'':=s'$. 
\end{proof}

Let $f: B\rightarrow P$ and $f': B'\rightarrow P$ be two maps. Set $f'\leq_{in} f$ if $B'\leq_{in} B$ and $f'(s')=f(s)$ for every $s' \in B'$, $s\in B$ with $s'\leq_{in} s$. Let   $\mathcal{H}_X(P)$ be the set  of maps $f: B\rightarrow P$  for which $B$ is a thin block with domain  $X$.
 
 \begin{fact} \label {fact2}
 Let $f: B\rightarrow P$ and $f': B'\rightarrow P$ with $f'\leq_{in}f$. If  $B' $ and $B$ are thin blocks then  $B'$ extends to a thin block $B''$ and $f'$  to a map $f''$ such that $\bigcup B''=\bigcup B$ and $f''\leq _{in} f$.
  \end{fact}
  \begin{proof} Applying $(e)$ of Fact \ref{fact1}, set $B'':= B'\cup D$ and define $f''$  by setting $f''(s''):=f(s)$ if $s''\in D$ and $f''(s''):=f'(s'')$ if $s''\in B''$. Then $f''\leq_{in} f$.
  \end{proof}
  
  \begin{fact}\label{fact3}   
   Let $P$ be a poset and $X\in \N^\omega$,  then:
  \begin{enumerate} [{(a)}]
  \item The relation $\leq _{in}$ is an order on the collection of maps $f$ whose domain is a thin block and whose range is $P$. 
  \item Every $\leq_{in}$-chain has an infimum on the set $\mathcal{H}_X(P)$.
  \item An element $f$  is minimal in $\mathcal{H}_X(P)$ if and only if every $f'$ with $f'\leq_{in} f$  is the restriction of $f$ to a  sub-block of the domain of $f$.
  \item If  $f$ is minimal in $\mathcal{H}_X(P)$ and   $f'\leq_{in} f$ has domain $C$ then $f^\prime$  is minimal  in $\mathcal{H}_{\bigcup{C}}(P)$.
  \item Let $B$ be a thin block  and $f:B\to P$. If $f$ is bad and $f'\leq_{in} f$ then $f'$ is bad.
  \end{enumerate}
  \end{fact}

\begin{proof}
$(a)$ Obvious.

$(b)$ Let $\mathcal D:= \{f_{\alpha}: B_{\alpha} \rightarrow P, \alpha\}$ be a chain of maps. Let $\mathcal C:= \{dom (f): f\in \mathcal D\}$. Then $\{s \in  \bigcup \mathcal C: s'\in  \bigcup \mathcal C \Rightarrow s' \not <s\}$ is a thin block  and  the infimum of $\mathcal{C}$. For $s\in D$, let $f'(s)$ be the common value of all maps $f_{\alpha}$. This map is the infimum of $\mathcal D$.

$(c)$ Apply Fact \ref{fact2}.

$(d)$ Follows from $(c)$. 

$(e)$ Apply  $(d)$ of Fact \ref {fact1}. 
\end{proof}
 
 \begin{lemma} \label{minimalbad}  
 Let $f$ be a  map from a thin block $B$ into $P$ and let $\mathcal F:= \{f'\in \mathcal{H}_{\bigcup B} (P) : f'\leq_{in} f\}$.  Then  there is a minimal $f'\in \mathcal{F}$ such that $f'\leq_{in} f$.
 \end{lemma}
  \begin{proof} 
 Follows from  Fact \ref{fact3} $(b)$ using   Zorn's Lemma.  
 \end{proof} 
   
 \vskip 8pt

 Let $s,t\in [\N]^{s}$. Set $s\leq_{end}t$ if $\lambda(s)\leq \lambda(t)$ and $s_{*}=t_{*}$.
\begin{lemma}  \label{key}Let $f:B\rightarrow P$ a bad map. If $P$ is wqo and 
  $f$ is minimal then there is an end-closed  barrier $B'\subseteq B$ such that:
 \begin{equation}\label{eqendclosed2}
s<_{end}t \mbox {\; in\; } B' \Rightarrow f(s)<f(t) \mbox{\; in\;  } P. \end{equation}  \end{lemma}
    \begin{proof}
Let $B_1$ be  a an end-closed subbarrier of $B$ and let $C_1:= \{s\cdot(b) : s\in B_1, b\in \bigcup B_1, b>\lambda(s)\}$. Divide $C_1$ into three parts  $D_{i}$, $i<3$, with $D_{i}:=  \{ s'\cdot (ab)\in C_1 : f(s'\cdot (a))\rho_{i}f(s'\cdot(b))\}$ where $\rho_{0}$ is the equality relation, $\rho_{1}$ is the strict order $<$ and $\rho_{3}$ is $\not \leq $ the negation of the order relation on $P$. 

Since $C_1$ is a barrier,  
Nash-Williams 's partition theorem (Theorem \ref{partnash} $(b)$) asserts that  one of these parts contains a   barrier $D$. Let $X$ be an infinite subset of  $\bigcup C_1$ such that $D=C_1\cap [X]^{<\omega}$.

 The inclusion $D\subseteq D_{2}$ is impossible. Otherwise, let  $s\in B_{1}$ such that  $s\leq_{in} X$, set $Y:= X\setminus s_{*}$ and set $g(a):= f(s_{*}\cdot (a))$ for $a\in Y$. Then $g$ is a bad map from $X$ into $P$. This contradicts the fact that $P$ is wqo. 
 
 The inclusion $D\subseteq D_{0}$ is also impossible. Otherwise, set $B':= \{s': s'\cdot (a)\in B_1\cap [X]^{<\omega} $ for some $a\}$. For $s'\in B'$,  set $f'(s'):=f(s'\cdot(a))$ where $a\in X$. In this case $f'(s')$  is well-defined. Since $P$ is wqo and $f$ is bad,  the order type of $B_{1}$ is at least $\omega^2$, hence $B'$ is a barrier. 
 The map $f'$ satisfies  $f'\leq_{in}f$. According to Fact \ref{fact3} $(c)$, the minimality of $f$ implies that $f'$ is the restriction of $f$ to $B'$.  Since $B'$ is not included into $B$ this is it not the case. A contradiction.
 
 Thus we have $D\subseteq D_{1}$. Set $B':= B_1\cap [X]^{<\omega}$. Then (\ref{eqendclosed2}) holds.
 
  \end{proof}

\subsection {An application to $S_{\omega}(P)$}

We deduce Theorem \ref{nonprincipal} from the equivalence $(i)\Longleftrightarrow (ii)$ in the following result. Without  clause $(ii)$, the result is due to A.Marcone \cite{marcone2001}. Without Marcone's result our proof only shows that  under clause $(ii)$ $P$ is $\alpha$-bqo. This suffices to prove   Theorem \ref{nonprincipal} but the result below is more precise.
\begin{theorem}\label{nonprincipal2}
Let $\alpha$ be a denumerable ordinal and $P$ be a poset. Then the following properties are equivalent:
\begin{enumerate}[{(i)}]
\item  $P$ is $\alpha \omega$- bqo;
\item   $P$ is $\omega$-bqo and $S_{\omega} (P)$ is $\alpha$-bqo.
\item $\mathfrak{P}_{\leq \omega}(P)$ is $\alpha$-bqo
\item $\mathfrak{P}(P)$ is $\alpha$-bqo
\end{enumerate}
\end{theorem}
\begin{proof}
 $(i)\Rightarrow (iv)$. Let $B$ be a barrier with order type at most $\alpha$ and $f: B\rightarrow \mathfrak {P}(P)$. If $f$ is bad,  let  $f': B'\rightarrow P$ where $B':= B^2$ and $f'(s\cup t)\in f(s)\setminus {\downarrow\hskip -2pt \hskip -2pt f(t)}$. (See Equation \ref{downarrow}.) This map $f'$ is bad and the order type of $B'$ is at most $\alpha\omega$. 
 
$(iv)\Rightarrow (iii)$  Trivial.
 
$(iii)\Rightarrow (ii)$ $P$ and $S_{\omega} (P)$ identify to subsets of $\mathfrak{P}_{\leq \omega}(P)$,  hence are $\alpha$-bqo.
 
 $(ii)\Rightarrow (i)$ Induction on $\alpha$. 
 Suppose that $P$ is not $\alpha\omega$-bqo. Let   $\beta$ be the smallest ordinal such that $P$ is not $\beta $-bqo. Then $\beta\leq \alpha \omega$.
 
 {\bf Case 1}. $\beta= \alpha'\omega$. According to Marcone \cite {marcone2001} the implication $(iii)\Rightarrow (i)$ holds for all denumerable ordinals, hence there is a bad map $f':B'\rightarrow \mathfrak{P}_{\leq \omega}(P)$ for which $B'$ is a barrier of type at most $\alpha'$. Let  $X\in \mathfrak{P}_{\leq \omega}(P)$. Since $P$ is wqo, $\downarrow\hskip -2pt \hskip -2pt  X$ is a finite union of ideals according to  a theorem of   Erd\"{o}s and Tarski (1943), see \cite{fraissetr}. 
 
Hence there are a finite antichain $A_X$ and a finite set  $B_X$ of strictly increasing sequences such that $\downarrow\hskip -2pt \hskip -2pt  X=\downarrow\hskip -2pt \hskip -2pt  A_X \cup \downarrow\hskip -2pt \hskip -2pt  B_X$. Let $g: B'\rightarrow \mathfrak {P}_{<\omega}(P)\times \mathfrak{P}_{<\omega}(S_{\omega}(P))$ defined by $g(s')= (A_f'(s'), B'_f(s'))$. This map is bad. Hence, from $(b)$ of Lemma \ref{higbqo} there is a bad map from a subbarrier $B''$ of $B'$ into $P$ or into $S_{\omega}(P)$.  The latter case is impossible since  $S_{\omega}(P)$ is $\alpha$-bqo and so is  the former case according to the  the induction hypothesis.
  
  
  {\bf Case 2}.  Case 1 does not hold, that is $\beta= \omega^{\gamma}$  where  ${\gamma}$ is a limit ordinal it follows that $\beta\leq \alpha$. Let  $f: B\rightarrow P$ be a bad map where $B$ is a barrier of type $\beta$. 
  
According to 
Lemma~\ref{minimalbad} there is a minimal  $f': B'\rightarrow P$ with $\bigcup B'=\bigcup B$ and $f'\leq_{end} f$ and according to Fact \ref{fact3} $(e)$ the map $f'$ is bad. Since $P$ is wqo, Lemma \ref{key} applies. Thus $B'$ contains  a subbarrier $B''$ on which 
$s\leq_{end }t$ implies $f'(s)< f'(t)$.

Let $F: B''_{*} \rightarrow P$ be given by
$F(s'):= \downarrow\hskip -2pt \hskip -2pt  \{f'(t)\in P:   t\in B''$ and $s'\leq _{in} t \}$. 

{\bf Claim 1} $F(s')$ is a finite union of non-principal ideals of $P$. Since $P$ is wqo, every initial segment is a finite union of ideals.  Hence in order to show  that  $F(s')$ is a finite union  of non-principal ideals  it suffices to show that it contains no  maximal element. 
Let $x \in F(s')$. Let $t \in B''$ such that $s'\leq_{in} t$, $f''(t)=x$.  Let
$u\in B''$  such that $t<_{end}u$. Then $s'\leq_{in}u$ hence $f''(u)\in  F(s')$. From Lemma \ref{key} $f''(t)<f(u)$, proving our claim.

{\bf Claim 2 } $F$ is good. Indeed, 
since $S_{\omega}(P)$ is $\alpha$-bqo, it follows from $(e)$ of Lemma \ref{higbqo} that  the collection of finite unions of its members is $\alpha$-bqo.

Hence $f'$ is good. Indeed,  since $F$ is good, there are 
$s', t'\in B''_{*}$ such that $s' \lhd   t'$ and $F(s')\subseteq F(t')$. Let $a:=t'(l(s')-1)$ then $s'\leq_{in} s:= s'.(a)\in B''$ then $f'(s'.(a))\in F(s')$. Since $F(s')\subseteq F(t')$ there is some $t\in B''$ such that  $t'\leq_{in }t$ and $f'(s)\leq f'(t)$. Because $s\lhd t$ the map $f'$ is good. 

This contradicts the hypothesis that $f'$ is bad and finishes the proof of the theorem.

\end{proof}

\section{The set of ideals of a well-quasi-ordered-poset}\label{sec:ideals}

In this section, we illustrate the relevance of the notion of ideal w.r.t. well-quasi-ordering. 

Define a  topology on 
 $\mathfrak{P}(P)$.  A basis of open sets consists of  subsets  of  the form 
$O(F,G):=\{X\in \mathfrak{P}(P): F\subseteq X \text{  and  } G\cap X=\emptyset \}$, where $F, G$ are 
finite subsets of $P$. The topological closure  of $down(P)$ in $\mathfrak{P}(P)$ is   a Stone space which is homeomorphic to the Stone space of $Tailalg (P)$, the Boolean algebra generated by $up(P)$ . With the order of 
inclusion added the closure of $down(P)$, $\bigcup{down(P)}$,  is isomorphic to the Priestley space of 
$Taillat (P)$  \cite{bekkali-pouzet-zhani}.

Note that $\bold{I}_{ <\omega}(P)$ is the set of compact elements of $\bold{I}(P)$, hence $\mathcal{J}(\bold{I}_{ <\omega}(P))\cong \bold{I}(P)$. We also note that $\mathcal{J}(P)$  is the set of join-irreducible elements of $\bold{I}(P)$.

We have
\begin{lemma}\label{emptyset} $\emptyset \not \in \bigcup {down(P)}\Longleftrightarrow P\in
\bold{F}_{<\omega}(P)$.
\end{lemma}

\begin{lemma}\label{ideals1}
$down(P)\subseteq \mathcal {J}(P)\subseteq \bigcup {down (P)}\setminus\{\emptyset\}$. In
particular,   the topological closures in $\mathfrak{P}(P)$ of
$down(P)$  and $\mathcal {J}(P)$ are the same. 
\end{lemma}

A poset $P$ is {\it up-closed} if every  intersection of two 
members of $up(P)$ is a finite union (possibly empty) of members of $up(P)$.

\begin{proposition}\label{ideals}
The following properties for a poset $P$ are equivalent:
\begin{enumerate}[{(a)}] 
\item  $\mathcal {J}(P)\cup 
\{\emptyset \}$ is closed for the product topology;
\item  $\mathcal {J}(P)=\bigcup {down(P)}\setminus \{\emptyset\}$;
\item $P$ is up-closed;
\item $\bold {F}_{<\omega}(P)$ is a meet-semi-lattice; 
\item $Taillat(P)=\bold{F}_{<\omega}(P)\cup \{P\}$.
\end {enumerate}  
\end{proposition}

Let us recall that a topological space $X$ is {\it scattered} if every non-empty subset $Y$ of $X$ contains an isolated point with respect to the topology induced on $Y$. We have:
\begin{proposition} Let $P$ be a poset. If $P$ is well-quasi-ordered  then ${\mathcal J}(P)$ is  a compact scattered space whose  set of   isolated points coincides with $down(P)$. 
\end{proposition}
\begin{proof} 

{\bf Claim 1} $\mathcal {J}(P)=\bigcup {down(P)}$.

Indeed, since $P$ is wqo it is up-closed. Hence, from Proposition \ref{ideals}, $\mathcal {J}(P)=\bigcup 
{down(P)}\setminus \{\emptyset\}$. Again, since $P$ is wqo, $P\in 
\bold{F}_{<\omega}(P)$. Hence, from Lemma \ref {emptyset},  $\emptyset \not \in \bigcup {down(P)}$ proving that $\mathcal {J}(P)=\bigcup 
{down(P)}$, as claimed. 

{\bf Claim 2}  As a subspace of the Cantor space $\mathfrak {P}(P)$, $\bold I(P)$ are   compact and  scattered.

$\bold I(P)$ is  closed. To see that it is scattered, let 
 $X$ be a non-empty subset of  $\bold {I}(P)$.  Since $P$ is wqo, $\bold {I}(P)$ is well-founded. Select a minimal element $I$ in $X$. Let $G:= \min(P\setminus I)$ and $O(G, \emptyset):= \{I'\in \bold {I}(P): G\cap I'=\emptyset \} $. Since $P$ is wqo, $G$ is finite, hence  $O(\emptyset, G)$ is a clopen subset  of  $\mathfrak {P}(P)$. Since $O(\emptyset, G)\cap X= I$,  $I$ is isolated in $X$.
 
 {\bf Claim 3} Let $J\in \mathcal J(P)$, then $J$ is isolated in $\mathcal J(P)$ if and only $J$ is principal.
 
 Suppose that $J$ is isolated. Then there is a clopen set of the form $O(F,G)$ such that $O(F,G)\cap \mathcal J(P)=\{J\}$. Since $J$ is up-directed, there is some $z$ in $J$ which majorises $F$. Clearly, $\downarrow\hskip -2pt  z \in O(F,G)\cap \mathcal J(P)$, hence $J=\downarrow\hskip -2pt  z$, proving that $J$ is principal. Conversely, let $z\in P$. Let $G:=\min( P  \setminus \downarrow\hskip -2pt  z)$. Since $P$ is wqo, $G$ is finite. Hence $O(\{z\}, G)$ is a clopen set. It contains only $\downarrow\hskip -2pt  z$, proving that $\downarrow\hskip -2pt  z$ is isolated.
 \end{proof}

From this result,  $\mathcal J^{\neg\downarrow\hskip -2pt }:= \mathcal J(P)\setminus down(P)$,  the set of non-principal ideals of $P$, coincides with $ \mathcal J^1 {(P)}$, the first derivative of $\mathcal J(P)$ in the Cantor-Bendixson reduction procedure. Our main  result establishes a link between the   bqo characters of  $\mathcal J(P)$ and  $ \mathcal J^1 {(P)}$. This suggests to look at the  other derivatives.

\section{Minimal type posets }\label{sec:min}

Well-quasi-ordered  posets  with just one non-principal  ideal are easy to describe. Each is a  finite unions of ideals, all but one being finite.  The infinite one, called   a {\it minimal type} poset,  can be  characterized in several ways:
\begin{proposition}[\cite{pouzettr}]Let $P$ be an
infinite  poset. Then, the following properties are  equivalent:
\begin{enumerate}[{(i)}]
\item $P$ is wqo and all ideals distinct from $P$  are principal;
\item $P$ has no infinite antichain and all ideals distinct from $P$ are finite;
\item Every proper initial segment of $P$ is finite.
\item Every linear extension of $P$ has order type $\omega$.
\item $P$ is level-finite, of height  $\omega$, and 
 for each $n<\omega$ there is  $m<\omega$ such that each element of 
height at most $n$ is below every element of height at least $m$.
\item $P$ embeds none of the following posets: an infinite antichain; a chain of order type
$\omega ^{dual} $; a chain of order type $\omega +1$; the direct sum $\omega\oplus 1$ of a chain of
order type $\omega$ and a one element chain.
 \end{enumerate}
\end {proposition}
An easy way of obtaining posets with minimal type is given by the following corollary
\begin {corollary}
Let $n$ be an integer and $P$ be  a poset.  The order on $P$ is the intersection of $n$ linear
orders of  order type $\omega$  if and only if $P$ is the intersection of $n$ linear orders  and $P$ has minimal type.
\end {corollary}

Minimal type posets occur quite naturally in symbolic dynamic. Indeed, let $S: A^{\omega}\rightarrow A^{\omega}$ be the shift operator on the set $A^{\omega}$ of infinite sequences $s:= (s_{n})_{n<\omega}$ of members of a finite set $A$ (that is $S(s):= (s_{n+1})_{n<\omega}$).  A subset $F$ of $A^{\omega}$ is {\it invariant} if $S(F)\subseteq F$.  As it is well-known, every compact (non-empty) invariant subset  contains a minimal one. To a compact invariant subset $F$ we may associate the set $\mathcal A(F)$ of finite sequences $s:= (s_{0}, \dots, s_{n-1})$ such that $s$ is an initial segment  of some member of $F$. Looking as these sequences as words, we may order $\mathcal A(F)$ by the factor ordering: a sequence $s$ being a {\it factor} of a sequence $t$ if  $s$ can be obtained from $t$ by deleting an initial segment  and an end segment of $ t$. 

We have then 
\begin{theorem}
 $\mathcal A(F)$ has minimal type if and only if $F$ is a minimal compact invariant subset.
\end{theorem}

\section{Maximal antichains, "pred" and "succ"}\label{sec:pred}

Let $P$ be a poset. We consider both  $A(P)$ and $AM(P)$ to be  ordered by domination. The main result of this section will be that if $AM(P)$ $\alpha$-bqo  then $P$ is $\alpha$-bqo. 

Our first aim is to prove that if $AM(P)$ is well founded then $P$ is well founded. To this end we will associate with every element $x\in P$ an antichain $\varphi(x)$  and investigate the connection between $x$ and $\varphi(x)$. First the following:

 \begin{lemma}\label{lem:1} Let $P$ be a poset and $X\subseteq A(P)$. Then the following properties are equivalent:
 \begin{enumerate} [{(i)}]
 \item $X$ is the minimum of $AM(P)$.
 \item $X$ is a minimal element of $AM(P)$.
 \item $P=\uparrow\hskip-2pt\hskip -2pt X$.
\end{enumerate}
\end{lemma}
\begin{proof}
Implications $(iii)\Rightarrow (i)\Rightarrow (ii)$ are obvious. 

\noindent
$(ii)\Rightarrow (iii)$: Suppose for a contradiction that  $P\setminus \uparrow\hskip-2pt X\not=\emptyset$. For every element $z\in P\setminus \uparrow\hskip-2pt X$ there is an element $x\in X$ with $z<x$. Let $Y$ be a maximal antichain of the set $P\setminus \uparrow\hskip-2pt X$ and $Z=X\setminus \uparrow\hskip-2pt Y$. Then $Z\cup Y$ is a maximal  antichain which is strictly dominated by $X$. 

\end{proof}

Let $P$ be a poset. Let  $\mathcal{P}$ be the order that is  the set of pairs $\mathcal{P}:=\{(x,y) : x\leq y\}$.

For $S\subseteq P$ let $\min S:=\{x\in S : y\in S \text{ and } y\leq x \text{ implies } y=x\}$ and for    $x\in P$ let  $\Phi(x):=\{z\in P : z\not< x\}$ and let  $\varphi (x):= \min\Phi(x)$.   Note that $x\in \varphi(x)$.

\begin{lemma}\label{lem:2} $x<y \Longleftrightarrow \varphi(x)<\varphi (y) \text{ and }x\not\in \varphi(y)$ for all $x,y\in P$.
\end{lemma}
\begin{proof}
Suppose $x<y$, then $\Phi(y) \subset \Phi(x)$. Hence $\varphi (y):=\min \Phi(y)\leq_{dom} \min\Phi(x)=:\varphi(x) $. ($\leq_{dom}$ is the domination order.) Then $\varphi(x)<\varphi(y)$ because $\varphi(x)\ni x \notin \varphi(y)$.

Conversely, suppose $\varphi(x)<\varphi (y)$ and  $x\not \in \varphi (y)$. If $x\not < y$, then  the definition of $\varphi (y)$ insures that $x'\leq x$ for some $x'\in \varphi(y)$. Since  $\varphi(x)<\varphi (y)$, we have $x=x'$ proving $x\in \varphi(y)$,  a contradiction.

\end{proof}

Note that if $\varphi(x)$ is a maximal antichain for every $x\in P$ and $AM(P)$ is well founded  then we obtain, using Lemma \ref{lem:2}, that $P$ is well founded. Actually we will show below that if $AM(P)$ is well-founded then $\varphi(x)$ is a maximal antichain.

Let $x\in P$  and the set $\mathcal  F\subseteq \mathfrak{P}(P)$. Then $\mathcal F(x):= \{F\in \mathcal F: x\in F\}$ and   $Inc_{P}(x):= \{y\in P: x \text{ and } y \text{ are incomparable}\}$.  Note that $AM(Inc_{P}(x))\cong AM(P)(x)$. ($\cong$ is the order isomorphism between he maximal antichains of $Inc_{P}(x)$ and the maximal antichains of $AM(P)(x)$ both ordered under domination.)  
This implies that   if $AM(P)$ is well-founded then $AM(\Phi(x))$  is well-founded. Let $S$ be a minimal element of $AM(\Phi(x))$. It follows then from  Lemma \ref{lem:1} that $S$ is the minimum of $AM(\Phi(x))$ and $\uparrow\hskip-2pt S=\Phi(x)$. Hence $S=\phi(x)$.  That is:  $\varphi(x)$  is the least maximal antichain of $P$ containing $x$. 

Also, if $P$ is well-founded then $\Phi(x)$ is well-founded, hence $\uparrow\hskip-2pt \varphi(x)=\Phi(x)$ which in turned implies using Lemma \ref{lem:1} that   $\varphi(x)$ is the least maximal antichain of $P$ containing $x$.  Hence we established  the following Lemma:

\begin{lemma}\label{lem:3} If $AM(P)$ is well-founded or if $P$ is well founded  then for all $x,y\in P$:
\begin{enumerate} 
\item $\uparrow\hskip-2pt \varphi(x)=\Phi(x)$.
\item  $\varphi(x)$ is the minimum  of all  maximal antichains of $P$ containing $x$.

\item $x<y \Longleftrightarrow \varphi(x)<\varphi (y)$ and $x\not\in \varphi(y)$. 
\item $P$ is well founded.
\end {enumerate}
\end{lemma}

Associated with the quasi order  $(P;\leq_{pred})$ is the equivalence relation $\equiv$ equal to the set  $\{(x,y) : x\leq_{pred} y \text{ and } y\leq_{pred}x\}$. Let $(P;\leq_{pred})/\equiv $ be the quotient equipped with the order induced by $pred(P)$. Let $\pi$ be the canonical map of $(P;\leq_{pred})$ to $(P;\leq_{pred})/\equiv $.  For every subset $S$ of $P$ let $\overline{\pi}(S):=\{\downarrow\hskip-2pt p(s); s\in S\}$.

\begin{theorem} \label{lem:AMembed}
{\    }
\begin{enumerate}
\item The function $\overline{\pi}$ induces an embedding of $AM(P)$  into  $\mathbf{I}\bigl((P;\leq_{pred})/\equiv\bigr) $ and if $P$ is well founded then $(P;\leq_{pred})/\equiv $ embeds into $AM(P)$. 
\item $\overline{\pi}$ induces an embedding of $\mathcal J^{\neg \downarrow\hskip -2pt } (P)$ into   $\mathcal J^{\neg \downarrow\hskip -2pt }((P;\leq_{pred})/\equiv)$.
\item If $P$ has no infinite antichain then this embedding is surjective, hence $\mathcal J^{\neg \downarrow\hskip -2pt } (P)\cong \mathcal J^{\neg \downarrow\hskip -2pt }((P;\leq_{pred})/\equiv )$. 
\end{enumerate}
\end{theorem}

\begin{proof} 
Let $A,B\in AM(P)$ with $\overline{\pi}(A)\subseteq \overline{\pi}(B)$. For every $a\in A$ there is an element $b\in B$ so that $a$ and $b$ are related under $\leq$. Assume for a contradiction that $b<a$. Then $\pi(b)<\pi(a)$. Because $\overline{\pi}(A)\subseteq \overline{\pi}(B)$ there is a $c\in B$ with $\pi(a)<\pi(c)$. This implies, because $b<a$, that $b<c$ a contradiction. Hence $\overline{\pi}(A)\subseteq \overline{\pi}(B)$ implies that $A$ is less than or equal to $B$ in the domination order which in turn implies that if $\overline{\pi}(A)= \overline{\pi}(B)$ then $A=B$. If $A$ is less than or equal to $B$ in the domination order then $\overline{\pi}(A)\subseteq \overline{\pi}(B)$ and hence we conclude that $\overline{\pi}$ is an embedding of $AM(P)$ into $\mathbf{I}\bigl((P;\leq_{pred})/\equiv\bigr) $.

We have: $x\leq_{pred (P)}y$ if and only if $\Phi(x)\supseteq \Phi(y)$ if and only if $\varphi(x)\leq \varphi (y)$ in $AM(P)$. Hence  $x\leq_{pred} y \text{ and } y\leq_{pred}x$, that is $x\equiv y$, if and only if $\varphi(x)=\varphi(y)$. This establishes item 1.

In order to establish item 2 let $J,J^\prime\in \mathcal J^{\neg \downarrow\hskip -2pt } (P)$.

If $J\subseteq J'$ then clearly $\overline{\pi}(J)\subseteq \overline{\pi}(J')$.   The functions $\pi$ and  $\overline{\pi}$ are order-preserving. Hence $J\in \mathcal J(P)$ implies $\overline{\pi}(J)\in \mathcal J(Q)$. Since $\pi$ is strictly increasing $\overline{\pi}(J)\in \mathcal J^{\neg \downarrow }(Q)$ if and only if $J\in \mathcal J^{\neg \downarrow }(P)$.

Suppose $J\not \subseteq J'$. Let $x\in J\setminus J'$. Since $J$ is not principal, there is some $x'\in J$ such that $x<x'$. Since $J'$ is an intial segment, $x'\not \in  J'$. Assume for a contradiction that  $\pi(x') \in \overline{\pi}(J')$. Then  there is an $x''\in J^\prime$ such that $\pi(x')\leq \pi(x'')$ hence $x'\leq_{pred} x''$. Therefore   $x<x''$  follows from $x<x'$.  Since $J'$ is an initial segment, we have $x\in J'$ contradicting the choice of $x$.

This proves that $\pi(J)\not \subseteq \overline{\pi}(J')$. Hence $\overline{\pi}$  is an embedding. 

Item 3:  Let $K\in \mathcal J^{\neg \downarrow\hskip -2pt } (Q)$. Let $J:= \{ x\in P: \pi(x)\in K\}$.  The set $J$ is an initial segment of $P$ since $\pi$ is order preserving.  The set  $J$ is a finite union of ideals since $P$ has no infinite antichain; see Fact \ref{fact:union}.  Let   $J:=J_{1}\cup \cdots \cup J_{k}$. We have  $K= \overline{\pi}(J)= \overline{\pi}(J_1) \cup \cdots \cup \overline{\pi}(J_{k})$. From the fact that $K$ is an ideal it follows that $K= \overline{\pi}(J_i)$ for some $i$. Since $K$ is not principal, $J_i$  cannot be principal. 

\end{proof}

\vskip 8pt

We derive Theorem \ref{criticbqo} from the following result:
\begin{theorem}\label{thmsucc}
Let $P$ be a poset with no infinite antichain and $\alpha$ be a countable ordinal . 
The following properties are equivalent:
\begin{enumerate}  [{(i)}]
\item $P$ is $\alpha$-bqo.
\item $(P;\leq_{succ})$ is $\alpha$-bqo.
\item $(P;\leq_{pred})$ is $\alpha$-bqo.
\item $(P; \leq_{crit})$ $\alpha$-bqo. 
\item $AM(P)$ is $\alpha$-bqo.
\end{enumerate}
\end{theorem}

\begin{proof}
Implications $(i)\Longrightarrow (iv)\Longrightarrow (iii)$ follow from the sequence of inclusions
$(\leq)\,  \, \subseteq\, \,  (\leq_{crit})\, \subseteq (\leq_{pred})$ and the implication $(iv)\Longrightarrow (ii)$ follows from the inclusion $( \leq_{crit})\, \subseteq\, \, ( \leq_{succ(P)})$.

$(iii)\Longleftrightarrow (v)$ Suppose $AM(P)$ is $\alpha$-bqo then $P$ is well-founded according to Lemma \ref{lem:3} and hence  according to item 1 of  Theorem  \ref{lem:AMembed}, $(P;\leq_{pred})/\equiv$ embeds into $AM(P)$ implying that  $(P;\leq_{pred})$ is $\alpha$-bqo. Conversely, suppose $(P;\leq_{pred})$ is  $\alpha$-bqo. Then ${\bf I}_{<\omega}(P;\leq_{pred})/\equiv)$ is $\alpha$-bqo  according to item 3 of  Theorem \ref{lem:AMembed}. From item 1 of  Theorem \ref{lem:AMembed}, the poset  $AM(P)$ embeds into ${\bf I}_{<\omega}{<\omega}(P;\leq_{pred})/\equiv)$ and  hence is $\alpha$-bqo.

We prove  implications $(ii)\Longrightarrow (i)$ and $(iii)\Longrightarrow (i)$ in a similar way as we have proven the  implication $(ii)\Longrightarrow (i)$ in Theorem \ref{nonprincipal2}.

Induction on $\alpha$. 
Let $Q$ be equal to $(P;\leq_{succ})$ or equal to $ (P;\leq_{pred})$. Suppose that $Q$ is $\alpha$-bqo. Since $(\leq)\, \subseteq\, (\leq_{ succ}\, \cap\,  \leq_{pred})$, the partial order  $P$ is well-founded and since it has no infinite antichain it is wqo.  If $P$ is not $\alpha$-bqo there is a barrier 
$B$  of type  at most $\alpha$ and a bad map $f :B\rightarrow P$. From Lemma \ref{minimalbad} there is a minimal  $f': B'\rightarrow P$ such that $\bigcup B'=\bigcup B$ and $f'\leq_{end} f$. According to Fact \ref{fact3} $(e)$ the map $f'$ is bad. Since $P$ is wqo, Lemma \ref{key} applies. Thus $B'$ contains  a subbarrier $B''$ on which 
\begin{equation}\label{eq1} s<_{end }t \Longrightarrow f'(s)< f'(t).
\end{equation}  

Suppose $Q:=(P;\leq_{succ})$. Since $(P;\leq_{pred})$ is $\alpha$-bqo, $f'$ cannot be bad thus  there are 
$s, t  \in B''$ such that $s\lhd t$ and   $f'(s)\leq_{succ(P)} f'(t)$. Pick  $t'\in B''$ such that $t<_{end}t'$. From (\ref {eq1}) we have $f'(t)< f'(t')$ .  According to the definition of $(P;\leq_{succ})$, we have  $f'(s)\leq f'(t')$. Since   $s\lhd t'$ it follows that $f'$ is good for $P$. A contradiction.

Suppose $Q:= (P;\leq_{pred})$. For $s\in B''$ set  $s^{+}:= s_{*}.(a)$ where $a$ is the successor of $\lambda(s)$ in $\bigcup{ B''}$. Set $f'^{+}(s):= f'(s^{+})$.  Since $(P;\leq_{pred})$ is $\alpha$-bqo there are $s$ and 
$t$ such that $s\lhd t$ and $f'^{+}(s)\leq_{pred}f'^{+}(t)$, that is $f'(s^{+})\leq_{pred} f'(t^{+})$.
Since $s<_{end}s^{+}$ we have  $f'(s)<f'(s^{+})$. According to the definition of $(P;\leq_{pred})$, this gives  $f'(s) < f'(t^{+})$. Since $s\lhd t^{+}$, the function $f'$ cannot be bad. A contradiction.
\end{proof}


\section{Maximal antichains with a prescribed size}\label{maximal}

\subsection{Two  element maximal antichains}

\begin{definition} 
Let $P$ be a poset. The structure $(P(2);\leq)$ is defined on $P(2):=P\times 2$ so that:
\[
(x,i)\leq (y,j) \text{ if } 
\begin{cases}
 i=j     & \text{ and  $x\leq y$, or  }  \\
 i=0 \text{ and } j=1   & \text{ and there exist  incomparable elements }\\
 &\text{ $x^\prime, y^\prime\in P$  with $x\leq x^\prime$ and $y^\prime\leq y$   }.
\end{cases}
\]
\end{definition}
It is easy to see that $P(2)$ is a poset.
 
 \begin{lemma} \label{construct1}Every poset  $P$ embeds into the  poset  $AM_2(P(2))$. 
 \end{lemma}

 \begin{proof} 

\noindent
{\bf Claim 1.}  If $y\leq x$ then  $(x,0)$  and $(y,1)$ are incomparable in $P(2)$.  The converse holds if  for every $x <y$  there are two incomparable elements $x', y'$ such that $x\leq x'$ and $y'\leq y$. 

If $(x,0)$ and $(y,1)$ are comparable then necessarily $(x,0)<(y,1)$. In this case there are two incomparable elements $x',y'$ such that $x\leq x'$ and $y'\leq y$. But if $y\leq x$, we get $y'\leq x'$,  a contradiction.  Conversely,  suppose that $(x,0)$  and $(y,1)$ are incomparable. Then clearly, $x$ and $y$ are comparable. Necessarily, $x\leq y$. Otherwise $x<y$. But,  from the condition stated, we have $(x,0)\leq (y,1)$,  a contradiction. 

\vskip 4pt
For $x\in V$, set $X_{x}:=\{(x, 0),(x,1)\}$.

\noindent
{\bf Claim  2.}  $X_{x}\in AM_2(P(2))$.

The set $X_x$ is an antichain according to Claim 1. Moreover, every element $(x', i')$  different from $(x,0)$ and $(x,1)$ is  comparable to one of these two elements. Indeed, if  $x'$ is comparable to $x$, then $(x',i')$ is comparable to $(x,i')$. If $x'$ is incomparable to  $x$ then $(x',i')$ is comparable to $(x, \neg i')$ where $\neg i' \not =i'$. This proves that $X_x$ is maximal.

\noindent
{\bf Claim 3. } The map $x \rightarrow X_{x}$ is an embedding of $P$ into $AM_{2}(P(2))$. That is:
$$x\leq y \Longleftrightarrow X_{x}\leq  X_{y}.$$
Suppose $x\leq y$. 
 Then we have  $(x,0)\leq (y, 0)$ and $(x,1)\leq (y,1)$ proving $X_{x}\leq  X_{y}$. Conversely, suppose $X_{x}\leq  X_{y}$, that is $(x,0)\leq(y,i)$ and $(x,1)\leq (y, j)$ for some $i,j\in \{0,1\}$. Due to our ordering, we have $j=1$, hence $x\leq y$ as required.
 \end{proof}

  With this construction, a poset $P$ which is not $\alpha$-bqo but is $\beta$-bqo for every $\beta<\alpha$ leads to a poset $Q$ having the same property and for which neither $AM_2(Q)$ nor $\bigcup {AM_2(P)}$ is  $\alpha$-bqo. The simplest example of this situation is given below.

 \begin{theorem}\label{thm:222}
Let $P$ be a poset with no infinite antichain, and $\alpha$ be a denumerable ordinal,  then $AM_2(P)$ is $\alpha$-bqo if and only if $\bigcup {AM_2(P)}$ is $\alpha $-bqo. 
\end{theorem}
\begin{proof}
If $Q:=\overline {AM_2(P)}$ is $\alpha $-bqo then $A(Q)$ is $\alpha$-bqo. In  particular, $AM_2(Q)$ is $\alpha $-bqo. This set is simply $AM_2(P)$ and the conclusion follows.

For the converse, we prove a bit more. Let $\mathcal P$ be a subset of $[P]^{2}$. We quasi-order $\mathcal P$ as follows:
$X\leq Y$ if for every $ x\in X$ there is some $y\in Y$ such that $x\leq y$ and for every $y\in Y$ there is some $x\in X$ such that $x\leq y$.

Let $T$ be a subset of $P':= \cup \mathcal P$.

{\bf Claim}  If $\mathcal P$ is $\alpha$-bqo then   $T$ is $\alpha$-bqo.

 Let $f: B\rightarrow T$ be a map from a barrier $B$ of type at most $\alpha$ into $T$. For each $s\in B$, select a map $F(s): 2:=\{0,1\}\rightarrow P$ such that $f(s)\in rg(F(s))\in \mathcal P$.    For $s\in B$,   set  $p(s):=i$  if $ F(s)(i)=f(s)$ and  for  $(s,t)\in B\times B$, set $\rho_{(s,t)}:=\{(i,j)\in 2\times 2: F(s)(i)\leq F(t)(j)\}$. Note that since an order is transitive,  for $s,t, u \in B$ the composition of relations satisfies

  \begin{equation}\label{eqcomp}
\rho_{(t ,u)}\circ \rho_{(s ,t)} \subseteq  \rho_{(s ,u)}
\end{equation}

 {\bf Subclaim 1} We may suppose that:
\begin{enumerate}
\item $p(s)=i_0$ for all $s\in B$ and some $i_0\in 2$;
\item $\rho_{(s,t)}= \rho$ for all pairs $(s,t)\in B\times B$ such that $s\lhd t$ and some $\rho\subseteq 2\times 2$;
\item  for every $i\in 2$ there are some $j,j'\in 2$ such that $(i,j), (j',i)\in \rho$
\end{enumerate}
{\bf Proof of Subclaim 1}
Since the map $p$ takes only  two values, we get from the partition theorem of Nash-Williams $p$ is constant  on a subbarrier of $B$. With no loss of generality, we may suppose that this barrier is $B$ proving that $1$ holds. Similarly, the map which associate $\rho_{(s,t)}$ to each element $s\cup t\in B^2$ takes only finitely many values  hence, by the same token, this map   is constant on a subbarrier $C$ of $B^2$. Necessarily $C= B^{2}\cap [X]^{<\omega}$for some  $X\subseteq \overline B$.  For $B':=B\cap [X]^{<\omega}$  the condition stated  in $2$ holds. We may suppose $B'= B$. Finally, since $\mathcal P$ is $\alpha$-bqo, the map which associates $rg(F(s)$ to $s\in B$ cannot be bad. According to the partition theorem of Nash-Williams this map is perfect on a subbarrier . With no loss of generality, we may suppose this subbarrier equals to $B$. From this $3$ follows.
\endproof

We will prove that  $\rho$ is reflexive. With conditions $1$ and $2$ it follows that $f$ is perfect, proving our claim (indeed,  let  $s \lhd t$. From $1$,  $f(s)=F(s)(i_0)$ and $f(t)=F(t)(i_0)$,  from $2$  $\rho_{(s,t)}=\rho$. The reflexivity
of $\rho$ insures that $(i_0,i_0)\in \rho$ that is  $(i_0,i_0)\in\rho_{(s,t)}$ which amounts to $F(s)(i_0)\leq F(t)(i_0)$. This yields $f(s)\leq f(t)$ as required).

 If $\rho$ is not reflexive, then it follows from condition $3$ that $\{(0,1), (1,0)\}\subseteq \rho$. From now, on we will suppose this later condition fulfilled.

 We say that two elements $s_0, s_1 \in B$ are {\it intertwined} and we set $s_0\lhd_{\frac{1}{2}} s_1 $ if there is an infinite sequence   $X:= a_0< \dots a_n <\dots $ of elements of $\overline B$ such that $s_0<_{init} X_{even}$ and  $s_1<_{init} X_{odd}$, where $X_{even}:= a_0< \dots a_{2n} <\dots $ and $X_{odd}:= a_1< \dots a_{2n+1}<\dots$.  We set $B^{(\frac{1}{2})}:= \{(s_0,s_1):  s_0\lhd_{\frac{1}{2}} s_1 \}$ and  $B^{\frac{1}{2}}:= \{s_0\cup s_1: (s_0,s_1)\in B^{(\frac{1}{2})} \}$  where $s_0\cup  s_1$ denotes the sequence $w$ whose range is the union of the ranges of $s_0$ and $s_1$.

 We note that
 \begin{enumerate}
 \item if $w\in B^{\frac{1}{2}}$ then the pair $(s_0,s_1)\in B^{(\frac{1}{2})}$ such that $w= s_0\cup s_1$  is unique;
 \item  $B^{\frac{1}{2}}$ is a thin block;
 \item if $X:= a_0< \dots a_n <\dots $ is an infinite sequence of elements of $\overline B$, $Y:= _{*}X$, $s_0, s_1, s_2 \in B$ such that $s_0<_{init} X_{even}$, $s_1<_{init} X_{odd}$, $s_2<_{init}Y_{odd}$ then $s_0\lhd_{\frac{1}{2}} s_1 \lhd_{\frac{1}{2}} s_2$ and $s_0\lhd s_2$.

 \end{enumerate}

Let  $w:= s_0\cup s_1, w':= s'_0\cup s'1\in B^{\frac{1}{2}}$. We say that $w$ and $w'$  are equivalent if there is a map $g$  from  $\cup \{ rg(F(s_i)): i<2\}$ onto $\cup \{rg(F(s'_i)): i<2\}$ such that
 \begin{enumerate}
 \item $g\circ F(s_i)=F(s'_i)$ for $i<m$;
 \item $\rho_{(s_i,s_j)}=\rho_{(s'_i,s'_j)}$ for all $i, j<2$;
  \end{enumerate}
 As one can check easily, this is an equivalence relation on $B^{\frac{1}{2}}$. Furthermore, the number of equivalence classes is finite (one can code each equivalence classe by a relational structure on a set of at most $4$ elements, this structure been made of four binary relations  and four unary relations). Since $B^{\frac{1}{2}}$ is a  block, it follows from the partition theorem of Nash-Williams that one class contains a barrier. Let $C$ be such a barrier,  $X\subseteq \overline B$ such that $C:= B^{\frac{1}{2}}\cap [X]^{<\omega}$ and  let $B':=B\cap [X]^{<\omega}$. For  $ s_0, s_1\in B'$ such that $s_0\lhd_{\frac{1}{2}} s_1 $, $\rho_{ s_0, s_1}$ and   $\rho_{ s_1, s_0}$ are constant; let $\rho_{\frac{1}{2}}$ and $\rho'_{\frac{1}{2}}$ their common value.

 For the proof of the next subclaims,  we select  $s_0, s_1, s_2 \in B'$ such that $s_0\lhd_{\frac{1}{2}} s_1 \lhd_{\frac{1}{2}} s_2$ and $s_0\lhd s_2$; according to condition 3 above this is possible.

{\bf Subclaim 2} $\rho_{\frac{1}{2}}\circ \rho_{\frac{1}{2}}\subseteq \rho$

 {\bf Proof of Subclaim 2} We have $\rho_{\frac{1}{2}}=  \rho_{s_0,s_1}= \rho_{s_1,s_2}$ and  $\rho= \rho_{s_0,s_2}$. The claimed inclusion follows from the composition of relations. \endproof

 {\bf Subclaim 3} $\rho'_{\frac{1}{2}}=\emptyset$

 {\bf Proof of Subclaim 3} Suppose the contrary; let $(i,j)\in \rho'_{\frac{1}{2}}$.  Case 1.  $i=j$. Let $k\not = i$. We have $(k,i) \in \rho= \rho_{s_0,s_2}$ and $(i,i)\in\rho'_{\frac{1}{2}}= \rho_{s_2,s_1}= \rho_{s_1,s_0}$. By composing these relations, we get with \ref {eqcomp}  $(k,i)\in   \rho_{s_0,s_0}$ contradicting the fact that  $rg(F(s_0))$ is an antichain. Case 2. $i\not= j$.  then from $(i,j)\in\rho'_{\frac{1}{2}}= \rho_{s_2,s_1}= \rho_{s_1,s_0}$ and $(j,i)\in \rho= \rho_{s_0,s_2}$ we get, by composing these relations, $(i,j)\in   \rho_{s_1,s_1}$  contradicting the fact that  $rg(F(s_1))$ is an antichain and proving Subclaim 3.\endproof

 {\bf Subclaim 4} $\rho_{\frac{1}{2}}$ satisfies condition $3$ of Subclaim $1$.

 {\bf Proof of Subclaim 4} Since $rg(F(s_0))$ and $rg(F(s_1))$ are two maximal antichains, each element of one is comparable to some element of the other. Since  $\rho_{s_1, s_0}=\rho'_{\frac {1}{2}}=\emptyset$, $rg(F(s_0))\leq rg(F(s_1))$ and the result follows.
 \endproof.

 Now,  if  $\rho_{\frac{1}{2}}$ is reflexive, it follows from Subclaim 2 that $\rho$ is reflexive and our claim is proved. If  $\rho_{\frac{1}{2}}$ is not reflexive then from Subclaim 4 it follows that $\{(0,1), (1,0)\}\subseteq \rho_{\frac{1}{2}}$. With Subclaim 4 this yields $(0,0), (1,1)\in \rho$ that is $\rho$ is reflexive and the proof of our claim is complete.

 \end{proof}

\subsection{ Rado's poset}

Let $V:=\{ (m,n) \in \N^2: m<n\}$.  We denote by $\leq_R$ the following relation on $V$:
\begin{equation}\label{eqrado}(m,n)\leq_R (m',n')\; \mbox {if either} \; m=m'\;  \mbox {and}\;  n \leq n'  \; \mbox {or}\;  n<m'
\end{equation}

This relation is an order. We denote  by $R$ the resulting poset. This poset, discovered by R. Rado \cite{rado}, is at the root of the discovery of bqo's.   R. Rado observed that $R$ is wqo but $I(R)$ is not wqo and has  shown that a poset $P$ is $\omega ^2$-bqo if and only if $I(P)$ is wqo.  R. Laver \cite {laver} has  shown that  a poset $P$ which is wqo, and not $\omega^2$-bqo  contains a copy of $R$. 
 Applying the  construction given in Lemma \ref{construct1} we have:
\begin{lemma}\label{counterexample}  
The poset $AM_{2}(R(2))$ is wqo but not $\omega^{2}$-bqo.
\end{lemma}
\begin{proof}
As a union of two wqo posets,  $R(2)$ is wqo. Hence $AM(R(2))$ is wqo for the domination order. In particular $AM_2(R(2))$ is  wqo.  Since $AM_2(R(2))$ embeds $R$, it cannot be $\omega^{2}$-bqo. 
\end{proof}

\begin{lemma} \label {AM(R)} 
$\bigcup {AM_{m}(R)}$ is bqo for every integer $m$ and $R$ embeds into $AM(R)$. 
\end{lemma}

\begin{proof} {\  $\, $}\\
\noindent
a) $\bigcup {AM_{m}(P)}$ is bqo. 
Let $m$, $m<\omega$.Then $\bigcup {AM_m(P)}\subseteq  \{(i,j): i<m, i<j<\omega\}$. Indeed,   let $A\in AM_{m}(P)$ then for each $i, i<m$ there is some $(i,j)\in A$ with $i<j$ (otherwise, add to $A$ an element $(i,j)$ with $j$ large enought). Consequently $\bigcup {AM_m(P)}$ is bqo. 

\noindent
b)  $R$ embeds into $AM(R)$. Since $R$ is  not $\omega^2$ bqo, $AM(R)$ is not $\omega^2$ bqo (Theorem \ref{thmsucc} ). Hence from Laver's  result mentioned above, the poset $R$ embeds into $AM(R)$. For the sake of simplicity we give a direct proof. 

Set $X_{(0,1)}:=\{(0,1),(1,2)\}$, $X_{(0,n)}:=\{(m,n): m<n\}$ for $n\geq 2$, $X_{(m,n)}:=\{(m',m): m'<m\}\cup \{(m,n)\}$ for $m\geq 1$.
One has to  check successively that: 

{\bf Claim 1.} $X_{(m,n)}$ is the least antichain in $AM(R)$ which contains $(m,n )$. 

{\bf Claim 2.} $(m,n)\leq (m',n')\Rightarrow X_{(m,n)}\leq X_{(m',n')}$.

{\bf Claim 3} $m,m'\geq 1$ and $X_{(m,n)}\leq X_{(m',n')}$ imply  $(m,n)\leq (m',n')$. 

\end{proof}

\subsection{Three element maximal antichains}

\begin{lemma}\label{construct2} Let $P:=(V;\leq)$ be  a poset. Let $L:=(V;\sqsubseteq)$ be a linear extension of $P$ with the property that if $x<y$ then there is a $z$ with $x\sqsubset z\sqsubset y$ and $z$ is incomparably in $P$ to both $x$ and $y$.    Then there is a poset $Q$ which is a union of a copy of $P$ and two copies of $L$ for which $\bigcup {AM_3(Q)}=Q$ and $AM_3(Q)$ is isomorphic to $L$. 
\end{lemma}

\begin{proof} 
  On $V\times 3$ define the following strict order  relation $<_Q$:
\[
(x,i)<_Q(y,j) \text{ if }
\begin{cases}
i=j=1     & \text{ and    $x < y$,  or}  \\
1\not=i \text{ or } j\not=1   & \text{ and $i\leq j$  and  $x\sqsubset y$. } 
\end{cases}
\]
Let  $Q=(V\times 3;\leq_Q)$ be the resulting poset by adding the identity relation to $<_Q$.   The order induced by $\leq_Q$  on   $V\times \{i\}$ coincides with the  order  $\leq$ on $V$ if $i=1$, whereas it coincides with $\sqsubseteq$ if $i\not = 1$.

Let $A:=\{(x_0,i_0), (x_1,i_1), (x_2,i_2), \dots, (x_{n-1},i_{n-1})\}$ be a finite antichain of $Q$ with $x_0\sqsubseteq x_1\sqsubseteq x_2\sqsubseteq x_3\sqsubseteq \dots\sqsubseteq x_n$.  If $i_j\not=2$ for any $j\in n$ then $\{(x_0,2)\}\cup A$ is an antichain of $Q$.  If $i_j\not=0$ for any $j\in n$ then $\{(x_{n-1},0)\}\cup A$ is an antichain of $Q$.  It follows that every element of $AM_3(Q)$ is of the form $\{(x_0,0), (x_1,1), (x_2,2)\}$ with $x_0\sqsupseteq  x_1\sqsupseteq  x_2$.

Let $A:=\{(x_0,0), (x_1,1), (x_2,2)\}\in AM_3(Q)$. Assume for a contradiction that $x_0\not= x_1$. Because $A$ is a maximal antichain it follows that $x_1<x_0$. According to the assumptions of the Lemma, there exists an element $y\in V$ which is not related to $x_1$ and $x_0$ and with $x_1\sqsubset y \sqsubset x_0$. Then $\{(x_0,0), (y,1), (x_1,1), (x_2,2)\}$ is an antichain. In a similar way we obtain that $x_2=x_1$. It follows that $AM_3(Q)=\{\{y\}\times 3 : y\in V\}$. 

We conclude that $AM_3(Q)$ is isomorphic to $L$ and $\bigcup AM_3(Q)$=Q.

\end{proof}

\begin{corollary}\label{construct2}
There exists a poset $Q$ for which $AM_3(Q)$ is bqo but $\bigcup AM_3(Q)$ is not bqo.
\end{corollary}
\begin{proof} 
 A poset  $P$ which is wqo and not bqo but satisfies  the conditions of Lemma \ref{construct2} leads to a poset $Q$ which is wqo,  not bqo, and  for which $AM_3(Q)$ is bqo but $\bigcup {AM_3(Q)}$ is not bqo. One may take for $P$  Rado's example and for a linear extension $L$ the lexicographic order according to the second difference.    
\end{proof}

\section{Index, notation,  basic definitions and facts}\label{index}

Let $P$ denote a partially ordered set.

 \vskip 15pt 
$\underline{\text{ \large  Poset, qoset, chain, well-founded, wqo, well-quasi-ordered:}}$  
\vskip 5pt \noindent  
If $(P;\leq)$ is a partially ordered set, a {\em poset}, we will often just write $P$ for $(P;\leq)$. We write $a\leq b$ for $(a,b)\in\,  \leq$. A {\it qoset} is a quasi ordered set and a linearly ordered poset is a {\it chain}.  

A qoset $P$ is {\it well-founded} if it contains no
infinite descending chain 
$$
\cdots < x_n < \cdots < x_0
$$
and if in addition, $P$ contains no infinite antichain then it is {\it well-quasi-ordered} (wqo). If $P$ is a chain and well-quasi-ordered   then it is {\it  well-ordered}.

\vskip 10pt \noindent
$\underline{\text{ \large   Initial segment, principal,  ${\mathbf I}(P)$, $\bold{I}_{<\omega}(P)$,  $\downarrow\hskip -2pt  X$:  }}$
\vskip 5pt \noindent
A subset $I$ of $P$ is an {\it initial 
segment} (or is {\it closed downward}) if $x\leq 
y$ and $y\in I$ imply $x\in I$. We denote by ${\bf  I}(P)$ the set of initial segments of $P$ ordered by inclusion.   

Let $X$ be a subset of $P$, then:
\begin{align}\label{downarrow}
\downarrow\hskip -2pt  X:= \{y\in P: y\leq x \text{  for some  }  x\in X\}.
\end{align}
We say that $\downarrow\hskip-2pt X$ is generated by $X$. If $X$ contains only one element $x$,  we write $\downarrow\hskip -2pt  x$ instead of $\downarrow\hskip -2pt  \{x\}$.  An initial segment generated by a singleton is  {\it principal} and it is  {\it finitely generated} if it is generated by  a finite subset of $P$. We denote by $\bold{I}_{<\omega}(P)$ the set of  finitely generated initial segments.

\vskip 10pt \noindent
$\underline{\text{ \large   $up(P)$, $down(P)$}}$
\vskip 5pt\noindent
We set  $up(P):=\{\uparrow\hskip-2pt\hskip -2pt  x: x\in P \}$ and 
$down(P):=\{\downarrow\hskip -2pt  x: x\in P\}$.

\vskip 10pt \noindent
$\underline{\text{ \large   $\leq_{dom}$, domination relation: }}$
\vskip 5pt \noindent
A subset $X$ of $P$ is being {\it dominated}  by the subset $Y$ of $P$, $X\leq_{dom}Y$,  if  for every $x\in X$ there is a $y\in Y$ such that $x\leq y$. The domination relation is  a quasi-order on the power-set $\mathfrak{P}(P)$. The  resulting ordered set is isomorphic to ${\bf I} (P)$, ordered by inclusion,  via the map which associates with $X\in \mathfrak{P}(P)$ the initial segment  $\downarrow\hskip-2pt X$.

\vskip 10pt \noindent
$\underline{\text{ \large   $S_{\omega}(P)$, strictly increasing sequence: }}$
\vskip 5pt \noindent
A sequence $(x_{n})_{n<\omega}$ of elements of  $P$ is {\it strictly increasing} if 
$$x_{0}<x_{1}\cdots <x_{n}<x_{n+1}<\cdots$$
We denote by $S_{\omega}(P)$ the set of strictly increasing sequences of elements of $P$. For  $(x_{n})_{n<\omega} \in S_{\omega}(P)$ and $(y_{n})_{n<\omega}\in S_{\omega}(P)$ we set $(x_{n})_{n<\omega} \leq (y_{n})_{n<\omega}$ if for every $n<\omega$ there is some $m<\omega$ such that $x_{n}\leq y_{m}$. This defines a quasi-order on $S_{\omega}(P)$. If we identify each $(x_{n})_{n<\omega}\in S_{\omega}(P)$ with  the subset $\{x_{n}:  n<\omega\}$ of $P$, this quasi-order  is induced by  the  domination relation on subsets.

\vskip 10pt \noindent
$\underline{\text{ \large    $\mathcal{J}(P)$, $\mathcal J^{\neg \downarrow\hskip -2pt }(P)$, ideal, non principal ideal:}}$
\vskip 5pt \noindent
An  {\it ideal} of  $P$  is a non empty initial segment  $I$ which is up-directed, that is every pair 
$x, y\in I$ has an upper bound $z\in I$.  Its {\it cofinality}, $cf(I)$, is the least cardinal $\kappa$ such that there is some set $X$ of size $\kappa$ such that  $I= \downarrow\hskip -2pt  X$.  The cofinality of $I$  is either $1$, in which case it has a largest element and  is said to be {\it principal}, or is infinite.  We denote by $\mathcal J^{\neg \downarrow\hskip -2pt }(P)$ the set of non principal ideals of $P$.

Note the following  fact, which goes back to Erd\"{o}s-Tarski (1943)(see \cite{fraissetr}):   
\begin{fact}\label{fact:union}
A  poset  $P$ has no infinite antichain if and only if every initial segment of $P$ is a finite union of ideals.
\end{fact}

\vskip 10pt \noindent
$\underline{\text{ \large    $P^{dual} $, $\bold {F}(P)$, $\bold{F}_{<\omega}(P)$, $\mathcal {F}(P)$, filter:}}$
\vskip 5pt \noindent
The {\it  dual} of  $P$ is the poset obtained from $P$ by reversing the order;  we denote it by $P^{dual} $. A subset
 which is respectively an initial segment, a finitely generated initial segment or an ideal of  $P^{dual} $ will
be called a {\it final segment}, a {\it  finitely generated final segment} or a {\it filter}
of
$P$. We denote  by
$\bold {F}(P)$, $\bold{F}_{<\omega}(P)$, and 
$\mathcal {F}(P)$ respectively,   
the collection of initial segments, finitely generated initial segments, and ideals of $P$  ordered by inclusion.

\vskip 10pt \noindent
$\underline{\text{ \large $\mathfrak{P}(P)$,  $A(P)$, $AM(P)$, $AM_n(P)$:}}$
\vskip 5pt \noindent
$\mathfrak{P}(P)$ denotes the set of subsets of $P$ and 
$A(P)$ is the collection of antichains of $P$, $AM(P)$ the collection of maximal antichains of $P$ and  $AM_n(P)$ is the collection of $n$-element maximal antichains of $P$. The  quasi-order of domination defined on $\mathfrak{P}(P)$ induces an ordering on the set $A(P)$ of antichains of $P$.
The sets $A(P)$ and $A(P^{dual})$ are equal. Hence, we may  order $A(P)$ by  the domination order of $P^{dual}$. In general, these two orders are distinct. But, they coincide on $AM(P)$.

\vskip 10pt \noindent
$\underline{\text{ \large   $\leq_{succ}$, $\leq_{pred}$, $\leq_{crit}$:}}$
\vskip 5pt \noindent
Let  $P$ be a poset. We write $a\leq_{pred(P)}b$  if $x<a$ implies $x<b$ for every $x\in P$.  We write  $a\leq_{succ(P)}b$ if $b<y$   implies $a<y$ for every  $y\in P \}$. We write $a\leq_{prec}b$, or $a\leq_{succ}b$, if $P$ is understood. We denote by $(P;\leq_{prec})$ and $(P;\leq_{succ})$ the corresponding quasi-ordered sets.

We write $a\leq_{crit(P)}b$  if $a\leq_{pred(P)}b$ and  $a\leq_{succ(P)}b$.  If in addition $a$ and $b$ are incomparable, the  pair 
$(a,b)$  {\it critical}.

\vskip 10pt \noindent
$\underline{\text{ \large   Interval order: }}$
\vskip 5pt \noindent
The poset  $P$ is an {\it interval-order} if  $P$  is isomorphic to a subset $\mathcal J$ of the set $Int (C)$ of non-empty intervals of some chain $C$. The  intervals are ordered as follows: for every $I, J\in Int(C)$,  $I<J$ if  $x<y$ for every $ x\in I$, every $y\in J$.  

Interval orders have  neat characterizations in different ways:  maximal antichains, associated preorders or  obstructions, see \cite {fishburn},  \cite{wiener}. We  recall this important characterization: 
\begin{theorem} The following properties are equivalent:
\begin{enumerate}[(i)]
\item  $P$ is an interval order.
\item  $Pred({ P})$ is total qoset.
\item  $Succ( P)$ is a total qoset.
\item  $P$ does not contain a subset isomorphic  to $\underline{2}\oplus \underline{2}$, 
the direct sum of two copies of the two-element chain.
\item 	$AM(P)$ is a chain.
\end{enumerate}
\end{theorem}

\vskip 10pt \noindent
$\underline{\text{ \large   $\N$, $[X]^{<\omega}$, $l(s)$, $\lambda (s)$, $_{*}s$, $s_{*}$, $s\cdot t$,  $s\leq_{in} t$,  $s\leq_{lex} t$:}}$
\vskip 5pt \noindent
The set of non-negative integers is denoted by $\N$,  the set of $n$-element subsets  of $X\subseteq \N$ by $[X]^n$  and the set of finite subsets of $X\subseteq \N$ by    $[X]^{<\omega}$.  We  identify each member $s$ of $[\N ]^{<\omega}$ with a  strictly increasing sequence, namely the   list of its elements written in an increasing order, eg $\{3,4,8\}$.  
Let $s\in [\N ]^{<\omega}$;
 the {\it length} of $s$,  $l(s)$, is    the number of its elements.  For  $m:= l(s)\not = 0$, we write $s:=\{s(0), \dots, s( m-1)\}$ with $s(0)<s(1)<\cdots<s(m-1)$.  The smallest element of $s$ is $s(0)$ the largest, denoted by $\lambda (s)$, is $s(m-1)$.

We denote  by $_{*}s$ the sequence obtained from $s$ by deleting its first element and by $s_{*}$ the sequence obtained by deleting the last element. (With the convention that $_{*}\emptyset=\emptyset_{*}=\emptyset$.) We denote by $(a)$ the one element sequence with entry $a$. 

Let  $s,t\in [\N]^{<\omega}$. If $\lambda(s) <t(0)$ then  $s\cdot t$ is the concatenation of $s$ and $t$. We denote by $s\leq_{in} t$ the fact that $s$ is an initial segment  of $t$ and by $s\leq_{lex} t$ the fact that $s$ is smaller than  $t$ in the lexicographic order. (For example $\{3,5, 8,9\}\leq_{lex} \{3,5, 9, 15\}$.) If there exists an $r\in [\N]^{<\omega}$ with $s<_{in} r$ and $t=_{*}r$ then $s\lhd t$.  For  example $\{ i\}\lhd\{j\}$ if and only if $i<j$, ($r=\{i,j\}$); also $\{i,j\}\lhd\{i',j'\}$ if and only if $j=i'$ and  then $r=\{i, i^\prime, j^\prime\}$.

\vskip 10pt \noindent
$\underline{\text{ \large    $\bigcup B$, $B_{\restriction X}$, block, thin block, barrier }}$
\vskip 5pt \noindent
 Let  $B\subseteq [\N ]^{<\omega}$ and $X\subseteq \N$.  Then   $\bigcup B:= \cup B$ and  $B_{\restriction X}:= B\cap [X]^{<\omega}$.  The set $B$ is a {\it block }\footnote{We stick to the  definition of Nash-Williams, 1968  \cite{nwbqots}; in some papers, a block is what we call a thin block} if: 
\begin{enumerate}
\item $B$ is infinite.
\item For every infinite subset $X\subset \bigcup B$ there is some $s\in B\setminus \{\emptyset\}$ such that $s\leq_{in} X$.
\end{enumerate}
If  $B$ is a block and  an antichain for the order $\leq_{in}$ then $B$ is a  {\it thin block}, whereas $B$ is a { \it barrier } if it is a block and  an antichain for the inclusion order.  A typical barrier is the set $[\N]^{n}$ of $n$-element subsets of $\N$. 

Trivially, every block contains a thin block,    the set $\min_{\leq _{in}}(B)$ of $\leq_{in}$  minimal elements of the block.  Moreover, if  $B$ is a  block, resp.a thin block,  and $X$ is an infinite subset of  $\bigcup B$ then $B_{\restriction X}$  is a block, resp a thin block.

\vskip 10pt \noindent
$\underline{\text{ \large   $_{*}B$, $B_{s}$, $_{s}B$, $B_{*}$, $B^2$, $B^{\circ}$, $B\leq_{in} B^\prime$:}}$
\vskip 5pt \noindent
Let $B, B^\prime\subseteq  [\N]^{\omega}$. 

$_{*}B:=\{_{*}s: s\in B\}$. For $s\in [ \N]^{<\omega}$ 

$B_{s}:= \{t \in  B: s\leq_{in} t\}$

$_{s}B:= \{t\in [\N]^{<\omega}: s\cdot t\in B\}=\{r\setminus s : r\in B_s\}$. 
Note  that if $B$ is a thin block and   $_{s}B$ is non-empty then it is a thin-block. 
 
$B_{*}:= \{s_{*}: s\in B\}$.  

$B^2:= \{ u:=s\cup t:  s, t \in B$ and $ s\lhd t\}$.   (This despite  the possible  confusion with the cartesian square of $B$.)

$B^{\circ}$ is the set of all elements $s\in B$ with the property that for all $i\in \bigcup{B}$ with $i<s(0)$ there is an element $t\in B$ with $(i)\cdot {_\ast s}\leq_{in} t$. In other words $s\in B^{\circ}$ if $(i)\cdot {_\ast s}\in \mathrm{T}(B)$ for all $i\in \bigcup{B}$ with $i<s(0)$. Let $B^\prime:=\{{_\ast s} : s\in B^\circ\}\setminus \{\emptyset\}$.
 
$B'\leq_{in} B$ if for every $s'\in B'$  there is some $s\in B$  such that  $ s'\leq_{in} s$. 
This is the quasi-order of domination associated with the order $\leq_{in}$ on $[\N]^{<\omega}$

\vskip 10pt \noindent
$\underline{\text{ \large   bqo, good, bad }}$
\vskip 5pt \noindent
A map $f$ from a barrier $B$ into a poset $P$ is {\em good} if there are $s,t\in B$ with $s\lhd t$ and  $f(s) \leq f(t)$. Otherwise $f$ is {\em bad}. \\
Let $\alpha$ be a denumerable ordinal. A poset $P$ is {\em $\alpha$-better-quasi-ordered}  if every map $f: B\rightarrow P$, where $B$ is a barrier of order type at most $\alpha$,  is good.  \\  
A poset $P$ is {\em better-quasi-ordered} if it is $\alpha$-better-quasi-ordered   for every denumerable ordinal $\alpha$.

\end{document}